\documentclass[12pt]{article}
\usepackage{graphicx}

\begin{document}

\title{Global action-angle variables for Duffing system}
\author{I.~Kunin$^1$ and A.~Runov$^2$}
\maketitle

\footnotetext[1]{
	 Department of Mechanical Engineering, 
	 University of Houston,
	 Houston, TX 77204, USA, e-mail kunin@uh.edu
	 }
\footnotetext[2]{
	 Department of Theoretical Physics, 
	 St.~Petersburg University,
	 Uljanovskaja~1, St.~Petersburg, Petrodvorez, 198504, Russia
	 }

\begin{abstract}
The classical representation of Hamiltonian systems in terms of action-angle
variables are defined for simply connected domains such as an interior of a
homoclinic orbit. On this basis methods of (local) perturbations leading, in
particular, to chaotic systems have been studied in literature. 

We are describing a new method for constructing global action-angle
variables and successive perturbations based on a topological covering of the
phase space. The method is demonstrated for representative example of the
Duffing system. 
\end{abstract}

The choice of variables for solutions of different problems is
typically related to matters of convenience. The question of ``the best''
variables looks from the first glance as not well posed. But starting from
the works of Poincare the special preference has been given to the variables
action (I) and angle ($\theta$) for Hamiltonian systems and their
perturbations~\cite{arn, golds}. At the same time the method of action-angle
variables
is typically restricted to simply connected domains, mainly to $R^n$. 

We are generalizing this method to global action-angle variables
defined globally for topologically nontrivial phase spaces. The approach is
based on topological transformations (covering) of the phase space plus
additional changes of geometry. 
In this publication the method is demonstrated for the popular conservative
and dissipative Duffing system. The (local) action-angle variables for the
system and their perturbations for chaos are considered, e.\,g.
in~\cite{meyer, wagg}.

\begin{figure}[p]
\centerline{
\includegraphics[width=10cm]{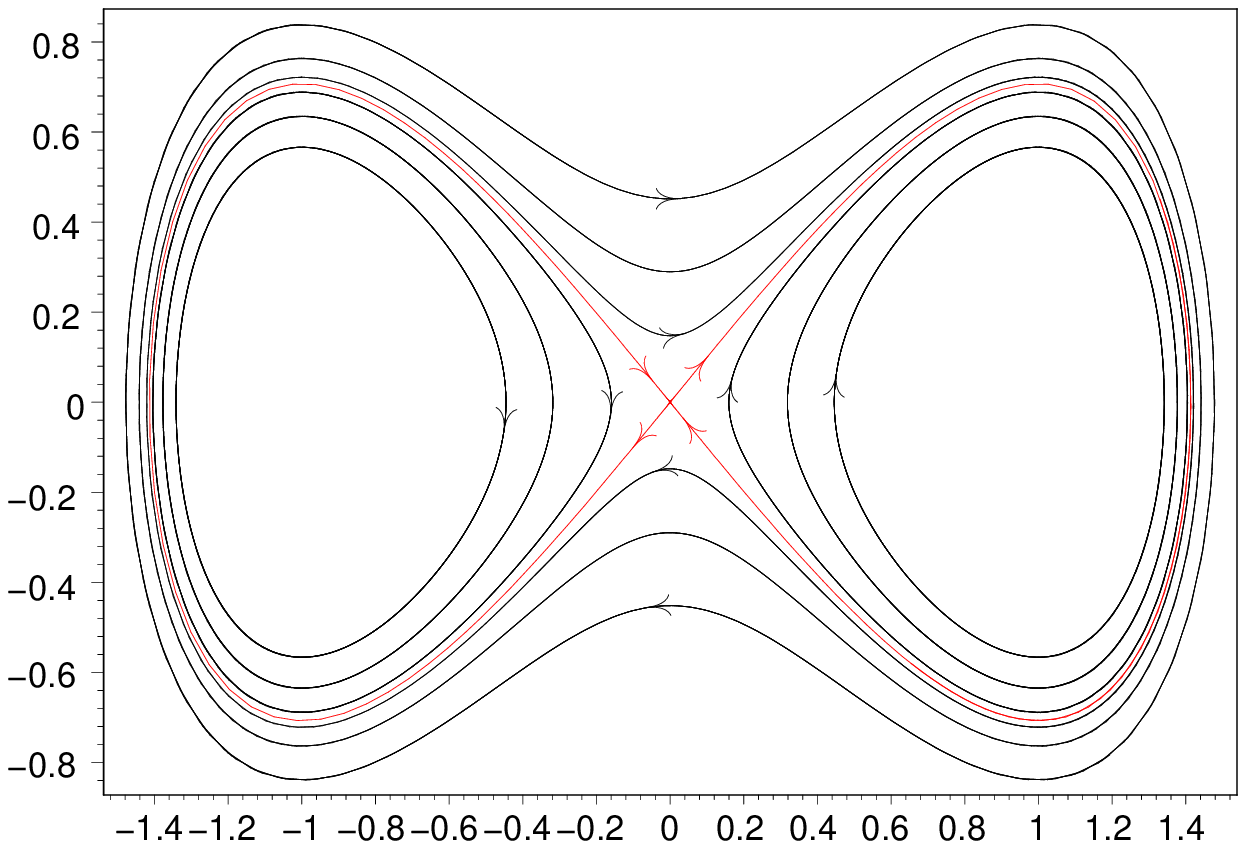}
}
\caption{Standard Duffing system.}
\label{stdDuff}
\end{figure}

Let us consider the well known Duffing equation (fig.~\ref{stdDuff}):
\begin{equation} \label{DE}
 \left\{
 \begin{array}{l}
  \dot x=y\\
  \dot y=x-x^3 -\mu y
 \end{array}
 \right.
\end{equation}
Let us begin with the case $\mu=0$. Then the corresponding Hamiltonian function is
\begin{equation}
 H(x,y)=\frac14x^4+\frac12y^2-\frac12x^2+C
\end{equation}
It is well known that one can introduce action-angle variables locally in each
of the three regions separated by the two homoclinic trajectories. Within the
conventional framework this cannot be done globally (for the entire plane)
because of the non-trivial topology of orbits. Let us utilize the symmetry of
the system with respect to the reflection $(x,y) \to (-x,-y)$. Due to this
symmetry, a half plane suffices for the description of the system. Consider
the following transformation, which maps a half plane into the plane:
\begin{equation} \label{coords}
 \begin{array}{l}
 x_1=x_1(x,y)=x^2-y^2,     \\
 y_1=y_1(x,y)=2xy.
 \end{array}
\end{equation}
Accordingly, the entire plane covers itself twice (a double covered plane).
With the help of a cut along the negative x-axis, we may speak of two planes
covering one. Let the points of the right half plane be mapped onto the upper
(red) plane and the points of the left half plane be mapped onto the lower
(green) plane. Then the formulas~(\ref{coords}) establish one-to-one correspondence
between the original plane $(x,y)$ and the double covered (colored) plane
$(x_1,y_1)$.

\begin{figure}[p]
\centerline{
\includegraphics[width=10cm]{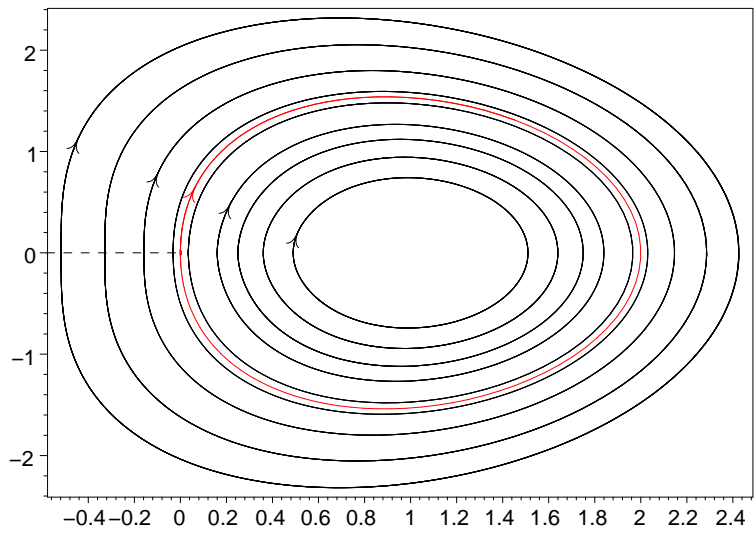}
}
\caption{Covered Duffing system.}
\label{covDuff}
\end{figure}

The key moment here is that the phase flow described by the equation~(\ref{stdDuff})  in
new variables does not depend on color, i.\,e. it is the same on upper and lower
planes.  This is an obvious consequence of the symmetry that was mentioned
above.  Correspondingly the equation~(\ref{stdDuff}) can be written in new variables in a
closed form
\begin{equation} \label{FDuff}
 \left \{
	 \begin{array}{l}
	 \dot x_1=\displaystyle\frac12(x_1+r) y_1	 \\
	 \dot y_1=\displaystyle -x_1^2-\frac12 y_1^2 +r(2-x_1)
	 \end{array}
 \right.
 , \qquad r=\sqrt{x^2+y^2},
\end{equation}
and information about covering is contained in the additional condition that
when the cut (negative $x$-axis) is crossed, the trajectories change over from
one plane to the other.   Because the phase flow is independent of color, this
transition is smooth (Fig.~\ref{covDuff}), and the motion assumes the character of
rotation around a single center.  Moreover, it is possible to consider a
projection of the system~(\ref{FDuff}) onto the $(x_1,y_1)$ plane which is
equivalent to gluing together two planes or washing out the colors. From a
more general point of view the resulting system may be interpreted as a
factorization of the original Duffing system by its symmetry group.

Returning back to the covered Duffing~(\ref{FDuff}) it is natural to rewrite
it in polar coordinates
\begin{equation}
 x_1=1+\rho \cos \theta, \qquad
 x_2=\rho \sin \theta,
\end{equation}
Then the angular velocity is
\begin{equation}
 \dot \theta=
-2\,{\frac {{x}^{6}+{x}^{4}{y}^{2}-2\,{x}^{4}+{y}^{4}+{x}^{2}+{y}^{2}}
{{x}^{4}+2\,{x}^{2}{y}^{2}+{y}^{4}-2\,{x}^{2}+2\,{y}^{2}+1}}
\end{equation}
It is strictly less than zero at all points except the origin where it is
equal to zero.  Accordingly, $\theta$ can be considered as new time
and explicitly
expressed through x and y:
\begin{equation}
\tan\theta(x,y)=y_1(x,y)/(x_1-1).
\end{equation}
In terms of $\theta$ and the corresponding action $I(H)$ the trajectories are
circles satisfying standard (but global!) action-angle equations.

\begin{figure}[p]
\centerline{
\includegraphics[width=10cm]{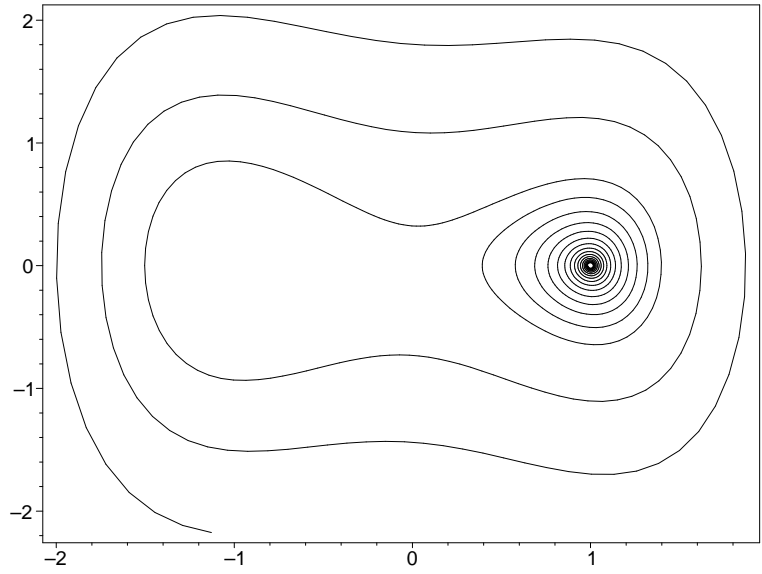}
}
\caption{Duffing system with dissipation, $\mu=0.1$}
\label{stdDuff_diss}
\end{figure}

\begin{figure}[p]
\centerline{
\includegraphics[width=10cm]{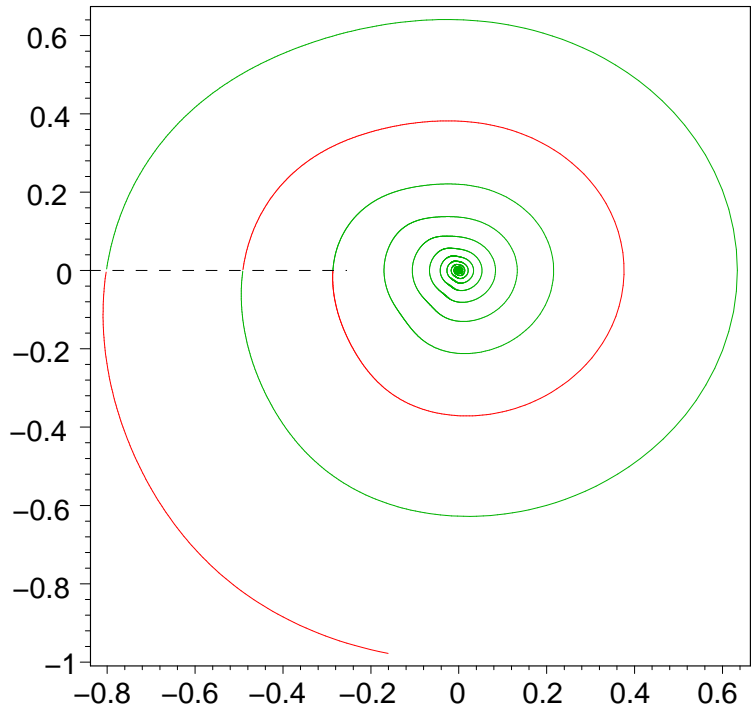}
}
\caption{Duffing with dissipation in $H$ --- $\theta$ variables}
\label{AA_diss}
\end{figure}

Now let us consider equation (\ref{DE}) with $\mu\ne0$ (Fig.~\ref{stdDuff_diss}).  
Then $\frac{dH}{d\theta}$ is not zero
any more, but is a strictly monotonous function:
\begin{equation}
 \frac{dH}{d\theta}=\frac{\dot H}{\dot \theta}>0,
\end{equation}
Accordingly, trajectories will be spirals (Fig.~\ref{AA_diss}).

This approach to global action-angle variables is in no way restricted to the
Duffing system. Its group theoretical background and a wider class of
applications will be considered elsewhere.

\end{document}